\documentclass[a4paper,10pt]{article}

\usepackage[english] {babel}
\usepackage{amsmath,amssymb,amscd}
\usepackage{amsthm}
\usepackage{verbatim}
\usepackage{authblk}

\usepackage[utf8]{inputenc}
\usepackage{ifthen}
\usepackage[all]{xy}

\newcommand{\N}{{\mathbb N}}

\newcommand{\R}{{\mathbb R}}

\def\set#1{\left\{\, #1 \,\right\}}
\def\abs #1{\left| \,#1\, \right|}
\def\norm #1{\left\| \,#1\, \right\|}
\def\norms #1{\| \,#1\, \|}
\def\inner #1#2{\langle \,#1,#2\, \rangle}
\def\EA0{E^A_{\,0}}

\def\calA{\mathcal{A}}

\def\calC{\mathcal{C}}

\def\calL{\mathcal{L}}

\def\calN{\mathcal{N}}

\def\calP{\mathcal{P}}

\theoremstyle{plain}
\newtheorem*{thmsn}{Theorem}
\newtheorem{thm}{Theorem}%[section]
\newtheorem{lemma}[thm]{Lemma}%[section]

\newtheorem{coro}[thm]{Corollary}
\newtheorem{prop}[thm]{Proposition}%[section]

\theoremstyle{definition}
\newtheorem{definition}{Definition}%[section]

\theoremstyle{remark}

\newtheorem*{notation}{Notation}

\def\noi{\noindent}

\begin{document}

\title{Geodesic rays of the $N$-body problem}

\author[$\dagger$]{J. M. Burgos}
\author[$\star$]{E. Maderna}

\affil[$\dagger$]{\small{CINVESTAV--CONACYT,
Ciudad de México, México.}}
\affil[$\star$]{IMERL, Universidad de la República, Uruguay.}

%\affil[$\dagger$]{
%Departamento de Matemáticas,
%Centro de Investigación y de Estudios Avanzados,\\
%Av. Instituto Politécnico Nacional 2508,
%Col. San Pedro Zacatenco,\\
%C.P. 07360 Ciudad de México, México}

%\email{burgos@math.cinvestav.mx}}
%\affil[$\star$]{
%Instituto de Matemática y Estadística Rafael Laguardia,
%Facultad de Ingeniería,\\
%Julio Herrera y Reissig 565,
%11300 Montevideo, Uruguay}
%\email{eze@fing.edu.uy}}

\maketitle

\begin{abstract}
For the Newtonian $N$-body problem,
we study the Jacobi-Maupertuis metric
of the nonnegative energy levels. 
We show that the geodesic rays are expansive,
that is to say, all the distances between the bodies
must be divergent functions. More precisely,
we prove that the evolution of such motions asymptotically decomposes
into free particles and subsystems in completely parabolic expansion.
The theorem applies in particular to the maximal characteristic curves
of any given global viscosity solution
of the stationary Hamilton-Jacobi equation $H(x,d_xu)=h$.

\medskip

\noi {\bf Keywords:} $N$-body problem, geodesic ray,
Jacobi-Maupertuis metric.

\noi {\bf AMS Subject Classification:}
70F10, 70F15, 37J50, 37B20, 70M20.
\end{abstract}
\bigskip

\section{Introduction}

The classification of the final evolutions
for solutions of the classical $N$-body problem
not experiencing singularities in finite time
(either collisions or pseudocollisions)
were developed during the 20th century,
starting with the pioneering works of Chazy \cite{Cha1,Cha2}.
The classification, in terms of the asymptotic behaviour
of the mutual distances between the bodies,
ends up being subordinated to a first description
that takes into account the order of growth of these distances.
More precisely, we are referring to the
Marchal-Saari primary description for the possible final evolution
of a Newtonian universe,
which applies to all motions, except those that may eventually
present an extremely complex behaviour.

Let us start by remembering this basic result. For this
we need to fix some notation. Let $E$ be an Euclidean space,
in which the $N$ bodies (punctual masses) evolve according
to the Newton's law of universal gravitation.
At every moment during the motion,
the different positions of the bodies form a configuration
$x(t)=(r_1(t),\dots,r_N(t))\in\Omega$,
where $\Omega\subset E^N$
denotes the open and dense set of configurations
\emph{without collisions}, that is,
such that $i\neq j$ implies $r_i\neq r_j$.
We will also denote by $R(t)=\max\set{\norm{r_i(t)-r_j(t)}\mid i<j}$
the greatest distance between two bodies at time $t$,
or equivalently,
the diameter of the subset of the Euclidean space $E$
defined by the configuration $x(t)$.
Then the Marchal-Saari theorem 
can be formulated as follows:

\begin{thmsn}[Marchal-Saari \cite{MarSaa}]
Let $x:[t_0,+\infty)\to E^N$ be a solution
of the $N$-body problem.
Either \[\limsup R(t)/t=+\infty\,,\]
or there is a configuration $a\in E^N$ such that
\[x(t)=ta+O(t^{2/3})\,.\]
\end{thmsn}

Motions in the first case of this alternative are
called \emph{superhyperbolic}.
To date, the existence of superhyperbolic motions
has not been proven, and it is known that their occurrence
require the presence of at least four bodies.
However, there is evidence that they could exist
and probably the elucidation of this dilemma will require an effort
of the same kind as that which recently concluded with
the resolution of the Painlevé conjecture by Xue \cite{Xue}.

In this paper, we address the study of
solutions to the $N$-body problem that,
in addition to being defined for all future time $t>t_0$,
also satisfy a strong global variational property,
namely of being geodesic rays.
As we will see, this variational property excludes
the possibility of superhyperbolicity,
and leads us naturally to the study of the $O(t^{2/3})$ term
established by the Marchal-Saari description.
To be more explicit,
we will now say what this variational property consists of.

As it is well known, negative energy levels are not projected
over the entire configuration space.
Indeed, the Newtonian potential is the function
\[U(x)=\sum_{i<j}m_im_j\,r_{ij}^{-1}\] where $m_i>0$ is the
value of the mass located at position $r_i\in E$,
$r_{ij}=\norm{r_i-r_j}$, and the kinetic energy of a
vector $v=(v_1,\dots,v_N)\in E^N$ is given by the quadratic form
\[K(v)=\frac{1}{2} \sum_{i=1}^N m_i\,\norm{v_i}^2\,.\]
Thus, the energy constant of a given motion $x(t)$ is
$h=K(\dot x(t))-U(x(t))$, which implies that the motion
must be confined to the \emph{Hill's region}
\[\Omega_h=\set{x\in\Omega\mid U(x)\geq -h}\,.\]
Clearly $\Omega_h=\Omega$ if and only if $h\geq 0$.
We also recall that trajectories of a fixed energy level $h$
are nothing other than the geodesics of the Jacobi-Maupertuis
Riemannian metric within the corresponding Hill's region.
More specifically,
this metric can be written as the conformal scaling
\[j_h=2(h+U)\,g_m\]
of the flat metric $g_m$ induced on $E^N$
by the mass scalar product.
The reader will find a full development of this geometric
point of view in the second author's recent work
with Venturelli \cite{MaVe2}.
In this article, the existence of geodesic rays is established
for arbitrary initial positions of the bodies
and for any choice of limit shape $a\in\Omega$,
thus giving rise to motions of the hyperbolic type.
Let us recall the notion of geodesic ray.
A \emph{geodesic ray} of a Riemannian manifold $(M,g)$
is an isometric embedding of
the half-line $[0,+\infty)$ within $M$, or in other words,
it is an arclength parameterized geodesic
$\gamma:[0,+\infty)\to M$ such that
all of its restrictions to compact subintervals
are minimizing geodesics.

Geodesic rays are essential elements for the study of
non-compact spaces,
and it is through them that Gromov introduced the
geometric notion of boundary at infinity.
However, taking into account that our motivation
is the study of dynamics,
we will give this notion a broader meaning.

\begin{definition}
\label{def:GR}
We will say that a curve $x:[0,+\infty)\to\Omega$
is a \emph{geodesic ray} for the metric $j_h$ whenever
\emph{its arclength parameterization} is an isometric embedding of
the half-line $[0,+\infty)$ into $(\Omega,j_h)$.
\end{definition}

In fact, we will see that this does not change things much,
since each geodesic ray has a unique arclength parameterization,
as well as it has a unique parameterization with constant energy $h$,
the latter being also defined over an unbounded time interval.
Moreover, this constant energy parameterization is indeed
a true motion of the $N$-body problem,
because it is a minimizer of the Lagrangian action functional
(see corollaries \ref{mainthmGR} and \ref{coro:ftmisgr} below).

We are now in position to summarize
the state of the art knowledge
on geodesic rays of the classical $N$-body problem,
as well as the contribution we will make in this paper.
We start by discussing their existence.
In what follows we assume, for all considered motions,
that the center of mass is at rest.

\begin{enumerate}
\item[(i)]
It is unknown wheather there are geodesic rays of negative energy levels.
It is not difficult to deduce that in every negative energy motion,
defined over an unbounded time interval,
at least one relative position $r_i(t)-r_j(t)$ between two bodies
must have a limit vector for a sequence of times $t_n\to +\infty$.
It seems reasonable to think that this form of recurrence could
contradict a global minimization property as strong as that
of being a geodesic ray.
However, we will not address this issue in this article.

\item[(ii)]
For $h=0$, the existence of geodesic rays
for arbitrary initial positions of the bodies
was first established in \cite{Ma1}, as calibrating curves
of any global viscosity solution of the Hamilton-Jacobi equation.
In \cite{PerSan}, Percino and Sánchez Morgado proved 
the existence of geodesic rays with prescribed
minimial central configuration as its limit shape.
On the other hand, it was proved in \cite{DaLMa}
that the geodesic rays of the zero energy level
are completely parabolic motions.
This means that $\dot x(t)\to 0$ as $t\to +\infty$,
or equivalently, that $r_{ij}(t)\approx t^{2/3}$ for all $i<j$.
As Chazy pointed out in \cite{Cha1},
the shape of the configuration of a completely parabolic
motion (as well as of a total collision motion)
must approximate the set of central configurations.
In fact, he introduced the notion of central configuration
for the general $N$-body problem
as the set of these limit shapes.
In \cite{Ma2},
the second author observes that
the set of central configurations being
limit shape of minimizing parabolic motions
also agree with the set of limit shapes of
minimizing total collision motions,
giving thus rise to the notion of
\emph{minimizing} central configuration.
It is not difficult to see that any minimal configuration
is minimizing.
However, it is unknown wheather there are
minimizing configurations that are not minimal.
Barutello and Secchi \cite{BarSec}
have established a sufficient condition to ensure
that certain configurations are not minimizing.
This result allows to prove that
for the three body problem with equal masses,
only the equilateral Lagrange configuration is minimizing,
and this holds for homogeneous potentials
in a certain range including the Newtonian case.
Moreover, the geodesic rays of the zero energy level of the planar
three body problem were studied by Moeckel, Montgomery
and Sánchez Morgado in \cite{MoMoSa}, where it is showed
for an open set of values of the masses,
that the only minimizing configuration
is the Lagrange equilateral one, and that
any completely parabolic motion with equilateral limit shape
admits a restriction which is a geodesic ray
(that is, the restriction to $t\geq t_0$ for some time $t_0$).

\item[(iii)]
For $h>0$ now it is known that there are, as we have said,
at least geodesic rays for arbitrary limit shape
\emph{without collisions} $a\in\Omega$ and starting
from any initial positions of the bodies \cite{MaVe2}.
These are therefore (completely) hyperbolic motions.
The existence of partially hyperbolic motions is part
of a current research by the first author, see \cite{Bur}.
These motions are those of the form $x(t)=ta+O(t^{2/3})$
such that their limit shapes have 	 collisions,
$a\in E^N\setminus\Omega$, but $a\neq 0$.
For this kind of motion there is associated
a natural cluster partition, defined as follows.
If $x(t)=(r_1(t),\dots,r_N(t))$ and
$a=(a_1\dots,a_N)$, then $a_i=a_j$ if and only if
$r_{ij}(t)=O(t^{2/3})$. The partition of the set of bodies
is then defined by this equivalence relation.
Therefore, the motion consists of a family of
clusters expanding between them with asymptotically
linear rate, while inside them the distances grow \emph{at most}
with a rate of order $t^{2/3}$.
For hyperbolic motions there are recent results by
Duignan, Moeckel, Montgomery and Yu \cite{DMMY},
in particular that the limit shape $a$ of a hyperbolic motion
is a real analytic function of the initial conditions of the motion.
\end{enumerate}

In this paper we prove that any geodesic ray of a positive
energy level is expansive, meaning that
$r_{ij}(t)\to +\infty$ for all $i<j$.
Clearly, we can also define these motions
as those for which $U(x(t))\to 0$.
Our main theorem will be formulated for free time minimizers of
the action functional
\[\calA_h(\gamma)=
\int_0^T\,L(\gamma(t),\dot\gamma(t))\,dt\;+\;hT\]
acting on the set of absolutely continuous curves
$\gamma:[0,T]\to E^N$.
We recall that as usual,
we say that a curve defined on a noncompact interval is minimizing
when all of its restrictions to compact subintervals are so.
Then, as a corollary, we will obtain that the same holds
for geodesic rays of the metric $j_h$.

\begin{thm}\label{mainthm}
Let $h>0$, and let $x:[0,+\infty)\to \Omega$
be a free time minimizer of the action functional $\calA_h$.
Then, writing $x(t)=(r_1(t),\dots,r_n(t))$ we have
\[\lim_{t\to +\infty}\norm{r_i(t)-r_j(t)} =+\infty\]
for all $1\leq i<j\leq N$.
\end{thm}

As we recall with more precision in the next sections,
the dynamics of the clusters, or subsystems,
was also analysed by Marchal and Saari in \cite{MarSaa}.
They prove that each subsystem behaves asymptotically
like a pure $P$-body problem $1<P<N$, in the sense that
energy and an angular momentum relationships
are asymptotically satisfied.

For geodesic rays of a positive energy level, we will prove that
each of the subsystems must be of zero asymptotic energy.
In particular, we will also be able to deduce,
using again the known results by Marchal and Saari,
that if one of the subsystem consists of $P$ bodies,
then the configuration which defines
must converge to the set of central configurations of
the $P$-body problem with the same value of the masses.

Actually, Theorem \ref{mainthm} is nothing but a trivial corollary
of the following one,
which gives a more accurate description
of the final evolution of geodesic rays of the
nonnegative energy levels.
In particular,
it generalizes the main theorem in \cite{DaLMa}.
The proof will be obtained essentially by
combining the decomposition of the Lagrangian action in terms of
the partial actions of the subsystems,
with the aforementioned results of Marchal and Saari
that we recall in section \ref{subsect:cluster-energy}.

\begin{notation}
In what follows we write $f\approx g$
whenever $f$ and $g$ are positive functions such that
their quotient is eventually bounded between two positive constants.
\end{notation}

\begin{thm}\label{mainthm0}
Let $h\geq 0$ and $x=(r_1,\dots,r_N):[0,+\infty)\to \Omega$ be
a free time minimizer of the action functional $\calA_h$.
Then for all $1\leq i<j \leq N$ we have
\[
r_{ij}(t)\approx t^{2/3}
\;\;\;\textit{or}\;\;\;
r_{ij}(t)\approx t\,.
\]
If all the distances grow like $t^{2/3}$ then $h=0$ and 
the motion is completely parabolic.
Otherwise $h>0$,
and the motion is completely hyperbolic
whenever all the distances grow like $t$.
\end{thm}

Roughly speaking, the theorem says
that every free time minimizer of $\calA_h$,
defined over an unbounded time interval,
can asymptotically be seen as a family of
subsystems which mutually separate each other
with a speed of higher order than that of
the separation of the bodies within them.
Moreover, it says that asymptotically,
each non trivial subsystem
(i.e. including at least two bodies)
can itself be considered as
a completely parabolic expanding system in uniform translation.

Now,
concerning the geometric interpretation let us say that,
for a curve defined over a compact interval,
it is well known the equivalence between the variational property
of being a free time minimizer of $\calA_h$,
and the property
of being a minimizing geodesic for the Jacobi-Maupertuis metric $j_h$.
However, we will see that the equivalence between,
on the one hand being a free time minimizer of $\calA_h$
defined over $[0,+\infty)$,
and on the other being a geodesic ray of the metric $j_h$
(parameterized with constant energy)
is in fact a consequence of Theorem \ref{mainthm}.
In this sense, we prove below (in Section \ref{sect: GR as FTM})
the following two corollaries.

\begin{coro}
\label{mainthmGR}
Let $h\geq 0$, and let $\gamma:[0,+\infty)\to\Omega$
be an arclength parameterized  geodesic ray
of the Jacobi-Maupertuis metric $j_h$.
Let $\gamma(s)=(q_1(s),\dots,q_N(s))$,
and $q_{ij}(s)=\norm{q_i(s)-q_j(s)}$.
Then we have:
\begin{enumerate}
\item[(i)] The parameterization of $\gamma$
with constant energy $h$ is a free time minimizer of $\calA_h$
and moreover, it is well defined for all $t\in [0,+\infty)$.

\item[(ii)]
If $h>0$, then for all $1\leq i<j \leq N$ either
\[
q_{ij}(s)\approx s^{2/3}
\;\;\;\textit{or}\;\;\;
q_{ij}(s)\approx s\,,
\]
and $q_{ij}(s)\approx s^2$ for all $1\leq i<j \leq N$ in the case $h=0$.
\end{enumerate}
\end{coro}

\begin{coro}
\label{coro:ftmisgr}
Let $h\geq 0$.
If $x:[0,+\infty)\to\Omega$ is a free time minimizer of
the action functional $\calA_h$,
then $x$ is a geodesic ray of the metric $j_h$.
\end{coro}

\section{The variational setting and preliminaries}

The preliminary aspects that follow will be approached
in a synthetic way.
For more precise details, or proofs of the statements,
we recommend that the reader consult mainly the references
\cite{DaLMa, MaVe2}, as well as those we will indicate.

As usual we define the Lagrangian as the function
$L:E^N\times E^N\to\R\cup \left\{+\infty\right\}$
\[
L(x,v)=K(v)+U(x)
=\frac{1}{2}\,\norm{v}^2+\sum_{i<j}\;\frac{m_im_j}{r_{ij}}
\]
where $\norm{v}^2=m_1\norm{v_1}^2+\dots+m_N\norm{v_N}^2$
for any $v=(v_1,\dots,v_N)$ and $r_{ij}$ is the distance
$\norm{r_i-r_j}$ between bodies $i$ and $j$
of the configuration $x=(r_1,\dots,r_N)$.
Of course $L<+\infty$ only for
$(x,v)\in T\Omega=\Omega\times E^N$.
It should be noted that the potential energy of the system
is the opposite of the Newtonian potential,
that is to say, $V=-U$.

The Lagrangian action
\[
\calA(\gamma)=\int L(\gamma(t),\dot\gamma(t))\,dt
\]
is well defined for
absolutely continuous curves in $E^N$ and it is clear that
for some curves $\gamma$ passing through collisions,
its value is $\calA(\gamma)=+\infty$.

It is well known, by experts in the field, that this does
not pose a problem for the application of variational methods,
but the fact that curves passing through collisions
can have \emph{finite} action.
This was already noticed by Poincaré,
claiming that the minimizers of the action
could a priori not be genuine solutions,
since they could present collision singularities.
Today there is no doubt that the resolution of this difficulty,
thanks to the celebrated Marchal's Theorem
\cite{Mar,Che,FerTer},
which assures that Poincaré's objection is not feasible,
has made it possible to fruitfully exploit variational methods.

To be explicit, we define for $x,y\in E^N$ and $\tau>0$
the sets of curves
\[\calC(x,y,\tau)=
\set{\gamma:[0,\tau]\to E^N
\textit{ absolutely continuous }\mid
\gamma(0)=x,\,\gamma(\tau)=y}\]
and
\[\calC(x,y)=\bigcup_{\tau>0} \calC(x,y,\tau)\,.\]
Usually in the literature, the actions are considered
on curves in Sobolev spaces.
It is not difficult to see that in our context,
this amounts to restricting our space of absolutely continuous
curves to that formed by those with finite action.
On the other hand, we prefer to keep the widest space,
since this is where Tonelli's theory naturally applies.
We also need to define, for each $h\in\R$,
the following \emph{action potentials}:
\[
\phi_h:E^N\times E^N\to \set{-\infty}\cup\R
\]
\[
\phi_h(x,y)=
\inf\set{\calA_h(\gamma)\mid \gamma\in\calC(x,y)},
\]
where $\calA_h(\gamma)$
is defined for a curve $\gamma\in\calC(x,y,\tau)$
as the action of $L+h$ on $\gamma$, that is to say,
\[
\calA_h(\gamma)=
\int_0^\tau\,L(\gamma(t),\dot\gamma(t))\,dt\;+\;h\tau\,.
\]
For $h<0$ it can easily be seen that $\phi_h=-\infty$.
For $h\geq 0$, It is well known that $(E^N,\phi_h)$
is a metric space, and indeed a length space.
Moreover, the restriction of $\phi_h$ to $\Omega$,
the space of configurations without collisions,
is precisely the Riemannian distance
induced by the Jacobi-Maupertuis metric $j_h$.

\subsection{Semistatic curves}
\label{subs: semistatics}

Before continuing, let us digress briefly
to show the relationship between
the motions we study here, and those that naturally arise in
the Aubry-Mather theory for Tonelli Lagrangians
(see for instance \cite{CoPa, Fat, Man}).
What is now clear,
is that for the $N$-body problem
the semistatic curves must escape,
hence both the Aubry and the Mather sets turn out to be empty.
However,
we find it is more suggestive to think that these sets
live in fact on the boundary at infinity for a convenient compactification,
and that eventually the inclusion of one in the other could be strict
for some particular values of the masses.

Following Mañé,
we define the critical value $c=c(L)$ of the Lagrangian as
\[
c=
\inf\set{h\in\R\mid \phi_h(x,y)>-\infty \text{ for all }x,y\in E^N},
\]

It follows that for $h<c$ we have $\phi_h=-\infty$
everywhere, and also that for any value of $h\geq c$,
the action potential $\phi_h$ defines a distance in $E^N$.
Moreover, for $h\geq c$ and for any pair of configurations
$x,y\in E^N$ there is always a curve
$\gamma\in\calC(x,y,\tau)$ for some $\tau>0$
realizing the $\phi_h$ distance between $x$ and $y$,
that is to say, such that $\calA_h(\gamma)=\phi_h(x,y)$.
As a consequence of Marchal's Theorem these curves
avoid collisions at any interior time and are true solutions.

A curve $\gamma$ defined on an interval $[t_0,+\infty)$
is called \emph{semistatic} whenever
\[
\calA_c(\gamma|_{[a,b]})=
\phi_c(\gamma(a),\gamma(b))
\]
for $c=c(L)$ the critical value,
and for any compact subinterval $[a,b]\subset [t_0,+\infty)$.
Observing that for any natural Lagrangian of the form $L=K-V$,
the critical value is nothing but $c=\sup V$,
we get in particular for the $N$-body problem that $c(L)=0$.
For this reason, here the \emph{critical action potential} is $\phi_0$,
and the semistatic curves are the free time minimizers
of the Lagrangian action.
Since these free time minimizers
are completely parabolic motions \cite{DaLMa} and
this type of motion have configurations whose shapes approximate
central configuration shapes,
their limit shapes define a very special class of central configurations.
Indeed, the set of \emph{minimizing} configurations
is defined in \cite{Ma2} precisely in this way.
On the other hand,
we have the set of \emph{minimal} configurations,
composed by those configurations which minimize
$U$ among all configurations with the same size.
Clearly any minimal configuration is a minimizing one.
This is because a minimal configuration is obviously
a central configuration,
and the corresponding parabolic homothetic motion
is a free time minimizer with fixed shape.  
This suggest that probably these two sets of configurations
will play an interesting role, respectively as
the Aubry and the Mather sets at infinity.

However, it should be said that for the positive energy levels
this analogy does not have such a clear interpretation .
Indeed for $h>0$, every configuration without collisions
$a\in\Omega$ can be obtained as the limit shape
of hyperbolic motions which are free time minimizers of $\calA_h$.

\subsection{Upper bounds for the action potentials}
\label{subs: estimates}

The following are estimates for minimal actions
between configurations.
All of them are deduced from the first one
which gives an upper bound for $\phi(x,y,\tau)$,
where $x$ and $y$ are given confgurations,
$\tau$ is a positive real number, and
\[
\phi(x,y,\tau)=
\min\set{
\calA_0(\gamma)=\calA(\gamma)
\mid\gamma\in\calC(x,y,\tau)
}.
\]
We will say that the configuration $x=(r_1,\dots,r_N)$
is \emph{contained} in the ball $B(r,R)\subset E$
whenever for all $i=1,\dots,N$ we have $\norm{r_i-r}<R$.

\begin{thmsn}[see \cite{Ma1}, Prop. 4, p. 1023]\label{mainestimate}
There are positive constants $\alpha_0$ and $\beta_0$ such that,
for any two configurations $x$ and $y$
contained in a given ball $B\subset E$ of radius $R>0$,
and for any $\tau>0$, there is a curve
$\gamma\in\calC(x,y,\tau)$ such that
\begin{enumerate}
\item[(i)]
for all $t\in[0,\tau]$
the configuration $\gamma(t)$ is contained
in the ball with the same center as $B$ but with radius $6NR$, and
\item[(ii)]
\[
\calA(\gamma)\leq
\;\alpha_0 \;\,MN^4\; \frac{\,R^2}{\tau}\;+\;
\beta_0 \;\,M^2N^4\; \frac{\,\tau}{R}\;
\]
where $M=m_1+\dots+m_N$ is the total mass of the system.
\end{enumerate}
\end{thmsn}

From this construction of canonical paths between configurations,
and using well chosen cluster partitions,
the following upper bounds in terms of $\norm{x-y}$ can be deduced.

\begin{thmsn}[see \cite{Ma1}, Prop. 9, p. 1031]\label{estim2}
There are positive constants $\alpha_1$ and $\beta_1$
which only depend on the number of bodies $N$ and the total mass $M$, 
such that for all $x,y\in E^N$ and any $\tau>0$ we have
\[
\phi(x,y,\tau)\leq
\;\alpha_1 \; \frac{\;\rho^2}{\tau}\;+\;
\beta_1 \;\frac{\,\tau}{\rho}\;
\]
whenever $\rho>\norm{x-y}$. 
In particular, in the case where $x\neq y$
the upper bound also applies for $\rho=\norm{x-y}$.
\end{thmsn}

Finally, from this last result we can obtain
the following upper bound for the action potentials
(see also Theorem 2 in \cite{Ma1} for the critical case $h=0$).

\begin{thmsn}[see \cite{MaVe2}, Theorem 2.11 p. 517]\label{estim3}
There are positive constants $\alpha$ and $\beta$
only depending on the number of bodies $N$ and the total mass $M$, 
such that for all $x,y\in E^N$ we have
\[
\phi_h(x,y)\leq
\left(\;\alpha\,\norm{x-y}\,+\,\beta h\,\norm{x-y}^2\;\right)^{1/\,2}
\]
for any value of $h\geq 0$.
\end{thmsn}

\section{Geodesic rays as free time minimizers}
\label{sect: GR as FTM}

So far we have assumed a total equivalence
up to reparameterizations of,
on the one hand the geodesic rays of the metric $j_h$
of a given energy level $h\geq 0$,
and on the other the free time minimizers
of the action functional $\calA_h$.
As we will see now,
this assumption deserves a non trivial justification.

For this, we will denote by $\calL_h(x)$ the length
of an absolutely continuous curve
$x:[a,b]\to\Omega$ with respect to the
$j_h$ metric, that is,
\[
\calL_h(x)=
\int_a^b \norm{\dot x(t)}_h\,dt=
\int_a^b [2(h+U(x(t)))]^\frac{1}{2}\norm{\dot x(t)}\,dt\,.
\]
Moreover, we recall that for any of these curves we have that
$\calA_h(x)\geq \calL_h(x)$,
with equality if and only if the curve $x$ has constant energy $h$.
 
\subsection{Proofs of corollaries \ref{mainthmGR} and \ref{coro:ftmisgr}}

\begin{proof}[Proof of Corollary \ref{mainthmGR}.]
Let $\gamma:[0,+\infty)\to\Omega$ be a geodesic ray
of $j_h$ for a given $h\geq 0$.
The reparameterization of $\gamma$ as a curve of constant energy $h$
is the curve $x:[0,t^*)\to\Omega$, with $t^*\in (0,+\infty]$,
defined by $x(t)=\gamma(\sigma(t))$ where
\[
\sigma:[0,t^*)\to [0,+\infty)
\]
is a diffeomorphism such that $\sigma(0)=0$, and such that
\[
\frac{1}{2}\norm{\frac{d\gamma}{ds}(\sigma(t))\,\dot\sigma(t)}^ 2
-U(\gamma(\sigma(t)))=h
\]
for all $t\in [0,t^*)$.
On the other hand, since we are assuming that $\gamma$
is arclength parameterized, for any $s\geq 0$ we have
\[
1=
\norm{\frac{d\gamma}{ds}(s)}^2_h=
j_h\left(\frac{d\gamma}{ds}(s),\frac{d\gamma}{ds}(s))\right)=
2(h + U(\gamma(s)))\norm{\frac{d\gamma}{ds}(s)}^2.
\]
Thus we deduce that $\sigma$ is well defined as the maximal solution
of the Cauchy problem
\[
\dot\sigma(t)=
2h + 2U(\gamma(\sigma(t)))
\]
with $\sigma(0)=0$, which clearly satisfies $\lim_{t\to t^{*}}\sigma(t)=+\infty$.

In order to prove the first statement, we observe that
for any $[a,b] \subset [0,t^*)$, we have
\[
\calA_h(x\mid_{[a,b]})=\calL_h(x\mid_{[a,b]})=
\calL_h(\gamma\mid_{[\sigma(a),\sigma(b)]})=
\sigma(b)-\sigma(a)
\]
because $x$ has constant energy $h$
and the length of a curve does not depend on the chosen parameterization.
Therefore, since $\gamma$ is a minimizing geodesic,
we conclude that $x\mid_{[a,b]}$ is a free time minimizer of $\calA_h$
for all $[a,b]\subset [0,t^*)$,
and it only remains to prove that $t^*=+\infty$.

Suppose $t^*<+\infty$.
Either $x$ has a limit configuration as $t\to t^*$,
or $x$ presents a noncollision singularity.
In both cases, a contradiction can be deduced
with the upper bound for the action potential $\phi_h$.
Since $\sigma(t)$ is precisely the $j_h$-length of $x\mid_{[0,t]}$,
if $x(t)\to x^*$ as $t\to t^*$, then for all $t\in [0,t^*)$ we have
\[
\sigma(t)=\phi_h(x_0,x(t)) < \phi_h(x_0,x^*)
\]
where $x_0=x(0)$,
contradicting the fact that $\sigma(t)\to +\infty$.
On the other hand, if we assume that $x$ presents a noncollision
singularity at $t=t^*$, then the Von Zeipel's Theorem \cite{McG,Zei}
implies that $\norm{x(t)}$ is unbounded on $[0,t^*)$.
Therefore $r_n=\norm{x(t_n)-x_0}\to +\infty$
for some sequence $t_n\to t^*$.
Let us denote by $A_n$ the Lagrangian action
of the restriction of $x$ to $[0,t_n]$, that is
\[
A_n=\calA(\,x\mid_{[0,t_n]}\,)=
\frac{1}{2}\int_0^{t_n} \,\norm{\dot x(t)}^2\,dt\;+\;
\int_0^{t_n} \,U(x(t))\,dt\;.
\]
By neglecting the potential term,
and applying the Cauchy-Schwarz inequality with the kinetic term,
we get the lower bound
\[r_n^2\;\leq
\left(\int_0^{t_n} \,\norm{\dot x(t)}\,dt\right)^2\leq
\;t_n\,\int_0^{t_n} \,\norm{\dot x(t)}^2\,dt \;<\;2t_nA_n\,,
\]
and therefore
\[
\frac{r_n^2}{2t_n}+ht_n < A_n+ht_n=\phi_h(x_0,x(t_n))
\]
for all $n\in\N$.
But the upper bound for the action potential cited in the previous section
(Theorem 2.11 in \cite{MaVe2}) gives
\[
\phi_h(x_0,x(t_n))\leq (\alpha\,r_n\,+\,\beta h \,r_n^2\,)^{1/\,2}
\]
and this inequality is incompatible with the previous one for $r_n$ big enough.
Hence we have proved that $t^*=+\infty$ as we wanted.

We prove now the second statement of the corollary.
We begin by considering the particular case $h=0$.
In this case, we have that $x(t)$ is a free time minimizer of $\calA_0$,
therefore, according to Theorem \ref{mainthm0},
all the mutual distances satisfy
\[
r_{ij}(t)\approx t^{2/3}\,.
\]
In turn, this implies that $\dot\sigma(t)=U(x(t))\approx t^{-2/3}$,
hence that $\sigma(t)\approx t^{1/3}$.
Now we define $\tau$
as the inverse function of $\sigma$,
and then we get that $\tau(s)\approx s^3$.
Since $q_{ij}(s)=r_{ij}(\tau(s))$,
we conclude that $q_{ij}(s)\approx s^2$ for all $i<j$.
On the other hand, if $h>0$ then we have
\[
\frac{\mu}{t}\leq U(x(t)) \leq \frac{\nu}{t^{2/3}}
\]
for some positive constants $\mu$ and $\nu$,
and for $t>0$ large enough.
Thus we get
\[
2ht +\mu\log(t) +\eta \leq \sigma (t) \leq 2ht + 3\nu\, t^{1/3} + \xi
\]
for $t$ large enough. We conclude that $\sigma(t)\approx t$
and therefore, its inverse function satisfies $\tau(s)\approx s$.
The proof is then achieved as before, by using Theorem \ref{mainthm0}
and the fact that $q_{ij}(s)=r_{ij}(\tau(s))$.
\end{proof}

\begin{proof}[Proof of Corollary \ref{coro:ftmisgr}]
Let $x:[0,+\infty)\to\Omega$ be a free time minimizer of $\calA_h$
and let $\gamma:[0,s^*)\to\Omega$
be the arclength parameterization of the curve,
where $s^*\in (0,+\infty]$. Thus we have
$\gamma(\sigma(t))=x(t)$ where
$\sigma:[0,+\infty)\to[0,s^*)$ is the length of $x\mid_{[0,t]}$.
Since $x$ is a free time minimizer we know that, for each $t>0$
\[
\phi_h(x_0,x(t))=\calA_h(x\mid_{[0,t]})=\calL_h(x\mid_{[0,t]})=\sigma(t)\,,
\]
which says that $\gamma$ is a minimizing geodesic of the metric $j_h$.
To see that $\gamma$ is a geodesic ray we only have to prove
that $\gamma$ has infinite length, that is $s^*=+\infty$.
Since a free time minimizer can not be superhyperbolic
(Lemma \ref{lemma: FTMnoSH} below),
we have that all mutual distances in the motion $x(t)$
satisfy $r_{ij}(t)=O(t)$.
Then we must have
\[
\frac{\lambda}{t}\leq U(x(t))
\]
for some constant $\lambda>0$ and $t$ large enough.
This suffices for $\sigma$ to be divergent, since
\[
\sigma(t)\geq \int_0^t U(x(\zeta))\,d\zeta\,.
\]
In particular $s^*=\lim_{\,t\to +\infty}\sigma(t)=+\infty$.
\end{proof}

\section{Systems and subsystems}

\subsection{The cluster's asymptotic energy}
\label{subsect:cluster-energy}

In the same quoted paper by Marchal and Saari,
the research continues with a finer study on the behaviour
of the subsystems that naturally appears in any non superhyperbolic motion.
In this section we recall some of these results in that direction
that we will use later in the proof of the main theorem.

Let us define first, for a motion $x(t)=ta+O(t^{2/3})$ as $t\to +\infty$,
its associated natural cluster partition,
as well as the geometry of that decomposition.

\begin{definition}
For a motion $x(t)=ta+O(t^{2/3})$ as $t\to +\infty$
its corresponding natural partition of the index set
$\calN =\set{1,\dots, N}$ is the one for which
$i,j\in\calN$ belong to the same class
if and only if the mutual distance $r_{ij}(t)$ grow as $O(t^{2/3})$.
Equivalently, if $a=(a_1,\dots,a_N)$, then
the natural partition is defined by the relation $i\sim j$ if and only if $a_i=a_j$.
The partition classes will be called clusters.
\end{definition}

\begin{definition}
Let $\calP$ be a given partition of $\calN$,
and let $x=(r_1,\dots, r_N)\in E^N$ be a given configuration.
For each cluster $A\in\calP$ we define the mass of the cluster
\[m_A=\sum_{i\in A}\,m_i\,,\]
the center of mass of the cluster
\[y_A=\frac{1}{m_A}\,\sum_{i\in A}\,m_ir_i\,,\]
and the relative positions $s_i=r_i-y_A$,
being $A$ the cluster containing $i\in\calN$.
\end{definition}

Notice that for $i,j\in A$ we have $s_{ij}=\norm{s_i-s_j}=r_{ij}$.

\begin{definition}
Let $\calP$ be any given partition of $\calN$.
Then for any given curve $\gamma(t)=(r_1(t),\dots,r_N(t))$
in $E^N$ we define, for each cluster $A\in\calP$, the functions
\[
K_A(t)=\frac{1}{2}\,\sum_{i\in A}m_i\norm{\dot s_i(t)}^2
\;\;\;\textrm{and}\;\;\;
U_A(t)=\sum_{i,j\,\in A,\,\;i<j}\;\frac{\,m_im_j\,}{s_{ij}(t)}
\;,
\]
as well as the function $H_A(t)=K_A(t)-U_A(t)$,
that we will call the energy of the cluster $A$ of the curve $\gamma$.
\end{definition}

\vspace{2cm}

\begin{thmsn}[Marchal-Saari \cite{MarSaa}, Theorem 2 p.165]
Let $x(t)$ be a motion defined for all $t\geq t_0$ which is not superhyperbolic.
If $\calP$ is its associated natural partition of $\calN$,
then for every cluster $A\in\calP$, there is a constant $h_A\leq 0$
such that
\[
H_A(t)=h_A +O(t^{-5/3})\,.
\]
\end{thmsn}

\begin{thmsn}[Marchal-Saari \cite{MarSaa}, Corollary 4 p.166]
Under the same hypothesis,
if $A\in\calP$ is such that $h_A=0$,
then $r_{ij}(t)\approx t^{2/3}$ for all $i,j \in A$ with $i\neq j$.
\end{thmsn}

Furthermore, we also have the following.

\begin{thmsn}[Marchal-Saari \cite{MarSaa}, Corollary 5 p.166]
Under the same hypothesis, if $h_A<0$ then the function
\[
r_A(t)=\min\set{r_{ij}(t)\mid i,j\in A,\;i\neq j}
\]
is bounded.
\end{thmsn}

\subsection{Geometry of the cluster partition}

Before proceeding with the proofs of the main results,
let us better describe the geometry of the cluster decomposition
of the configuration space associated with a given partition of the system.
We believe that a convenient notation is necessary
to be able to express this description synthetically.
For this purpose,
we start with the following definitions dependent 
on a given partition $\calP$ of the set $\calN =\set{1,\dots, N}$.

The space of the cluster centers will be the space
$E^\calP=\set{(u_A)_{A\in\calP}}$ endowed
with the inner product which weights each center
with the total mass of the cluster.
More precisely, for any pair $u,u'\in E^\calP$,
\[
\inner{u}{u'}=\sum_{A\in\calP}\;m_A\,\inner{u_A}{u'_A}_E
\]
where as before, $m_A=\sum_{i\in A}m_i$.
Second, for each cluster $A\in\calP$
we define the space of centered configurations of the cluster
\[
\EA0=\{\;z=(s_i)_{i\in A}\in E^A\;\mid \; \sum_{i\in A} m_i\,s_i=0\;\},
\]
also endowed with its corresponding mass inner product
\[
\inner{z}{z'}=\sum_{i\in A}\;m_i\,\inner{s_i}{s'_i}_E\,.
\]

Notice that for this product, $\EA0$ is precisely the orthogonal
complement of the diagonal $\Delta_A\subset E^A$.
Also notice that $\dim E^\calP=\#\calP\,\dim E$,
and that for each $A\in\calP$ we have $\dim \EA0=(\#A-1)\dim E$.

Finally,
we define the product of all these inner product spaces,
that is,
\[
D_\calP=\;E^\calP\,\times\;\prod_{A\in\calP}\,\EA0\,
\]
and actually, it is well known that this space is isometric
to the usual mass inner product space defined on the original
space of configurations $E^N$.

\begin{prop} For any partition $\calP$ of $\calN$,
the linear map
\[
d:E^N\to D_\calP
\]
given by
\[
x=(r_1,\dots,r_N)\;\mapsto\; d(x)=(y(x), (z_A)_{A\in\calP})\,,
\]
where $y(x)=(y_A)_{A\in\calP}$ is the configuration
produced by the centers of mass of the clusters,
and $z_A=(s_i)_{i\in A}$ are the configurations of
the clusters given by the relative positions $s_i=r_i-y_A$, is an isometry.
\end{prop}

\begin{proof}
It is trivial to see that the linear map $d$
takes values in the space $D_\calP$. This follows from the fact that
the center of mass of each configuration $z_A$ is precisely
\[
\sum_{i\in A}m_is_i=
\sum_{i\in A}m_i(r_i-y_A)=\sum_{i\in A}m_ir_i-(\sum_{i\in A}m_i)y_A=
m_Ay_A-m_Ay_A=0.
\]

Let us prove now that $d$ is an isomertic immersion.
Let $x,x'\in E^N$ be any two configurations.
If $d(x)=(y,(z_A)_{A\in\calP})$ and $d(x')=(y',(z'_A)_{A\in\calP})$
are the corresponding cluster decompositions, then we have

\begin{eqnarray*}
\inner{d(x)}{d(x')} & = &
\inner{y}{y'} +
\sum_{A\in\calP} \inner{z_A}{z'_A}\\
&=&
\sum_{A\in\calP}m_A\inner{y_A}{y'_A}+
\sum_{A\in\calP}\sum_{i\in A}m_i\inner{s_i}{s'_i}\\
&=&
\sum_{A\in\calP}\,\sum_{i\in A}m_i\inner{y_A}{y'_A}+
\sum_{A\in\calP}\sum_{i\in A}m_i\inner{r_i-y_A}{r'_i-y'_A}\\
&=&
\sum_{i=1}^Nm_i\inner{r_i}{r'_i}=\inner{x}{x'}\,.
\end{eqnarray*}

Finally, we conclude that $d$ is an isometry since
\[
\dim D_\calP = \#\calP\dim E+\sum_{A\in\calP}(\#A-1)\dim E=
N \dim E=\dim E^N\,.
\]
\end{proof}

In particular we deduce the following useful corollary.

\begin{coro}
For any partition $\calP$ of $\calN$, 
the decomposition
$d(x)=(y, (z_A)_{A\in\calP})$
of a given configuration $x\in E^N$
with respect to the cluster partition $\calP$ satisfies
\[
\norm{x}^2=\norm{y}^2\,+\,\sum_{A\in\calP}\,\norm{z_A}^2\,.
\]
\end{coro}

\subsection{Decomposition of the Lagrangian action}
\label{subsec:decomposition}

The cluster decomposition of the configurations will be used
to decompose the Lagrangian action of a curve.
For this, let us first observe that for any curve $\gamma(t)$ in $E^N$
with decomposition $d(\gamma(t))=(\xi(t),(\sigma_A(t))_{A\in\calP})$
we have $\dot\sigma_A(t)\in \EA0$
at any point of differentiability of $\gamma$.
In particular, for absolutely continuous curves the previous
lemma gives the almost everywhere identity
\[
\norm{\dot\gamma}^2=\norms{\dot\xi}^2\,+
\,\sum_{A\in\calP}\,\norm{\dot\sigma_A}^2\,.
\]

In order to decompose the Lagrangian action of a curve
we also need to decompose the Newtonian potential.
For that purpose, given a partition $\calP$ of $\calN$,
for any configuration $z\in\EA0$, say $z=(s_i)_{i\in A}$, 
we write
\[
U_A(z)=\sum_{i,j\,\in A,\,\;i<j}\;\frac{\,m_im_j\,}{s_{ij}}\,,
\]
where $s_{ij}=\norm{s_i-s_j}$. Morever, if $x=(r_1,\dots,r_N)\in E^N$,
for any pair of different clusters $A,B\in\calP$, 
we define the interaction potential $W_{AB}$ by
\[
W_{AB}(x)=
\sum_{i\in A,\,j\in B}\frac{\,m_im_j\,}{r_{ij}}\,.
\]
Now, if the cluster decomposition of $x$ is $(y_c, (z_A)_{A\in\calP})$,
 where $z_A=(s_i)_{i\in A}$,
then clearly we have $r_{ij}=s_{ij}$
for all $i,j\in A$ and any cluster $A\in\calP$.
Therefore we can write
\[
U(x)=
\sum_{A\in\calP}U_A(z_A)\;+\,W_\calP(x)
\]
where
\[
W_\calP(x)=\sum_{A,B\in\calP\;A\neq B}\;\frac{1}{2}\;W_{AB}(x)\,.
\]

Finally we conclude that for any absolutely continuous curve
$\gamma:[a,b]\to E^N$, the Lagrangian action decomposes as
\[
\calA(\gamma)=
\int_a^b\frac{1}{2}\norms{\dot\xi}^2\;+\;
\sum_{A\in\calP}
\;\left(\int_a^b\frac{1}{2}\norm{\dot\sigma_A}^2\,+\,U_A(\sigma_A)\right)\;+\;
\int_a^bW_\calP(\gamma)
\]
where $\xi(t)$ is the configuration of the cluster centers of $\gamma(t)$,
and $\sigma_A(t)$ is the centered configuration
defined by the cluster $A$ of $\gamma(t)$, for each $A\in\calP$.

\section{Proof of the main theorem}\label{sec: Proof}

We start by showing that curves in our hypothesis are not superhyperbolic.
This result was first established in the previous work
by the second author with Venturelli for the study
of calibrating curves of some viscosity solutions
of the Hamilton-Jacobi equation,
and we transcribe here the proof
for the sake of completeness.

\begin{lemma}[\cite{MaVe2}, p.534]\label{lemma: FTMnoSH}
If $x:[0,+\infty)\to E^N$ is a free time minimizers of $\calA_h$
with $h>0$, then $\gamma$ is not superhyperbolic.
\end{lemma}

\begin{proof}
Suppose we have a sequence $t_n\to +\infty$ such that
$\norm{x(t_n)}t_n^{-1}\to +\infty$.
We will prove that this is in contradiction with the upper bound
for $\phi_h$ given in Section \ref{estim3}.
Therefore we will conclude that $\norm{x(t)}<\lambda t$ for some
constant $\lambda$ and $t$ big enough, which says that
$x(t)$ is not superhyperbolic.

Let us call $r_n=\norm{x(t_n)-x(0)}$ and $A_n=\calA(\,x\mid_{[0,t_n]})$.
Then we have, as in the proof of Corollary \ref{mainthmGR},
on the one hand, that
\[
\frac{r_n^2}{2t_n}+ht_n \leq A_n+ ht_n = \phi_h(x(0),x(t_n))\,,
\]
and on the other hand,
because of the estimate given in section \ref{subs: estimates},
that
\[
\phi_h(x(0),x(t_n))\leq (\alpha\,r_n\,+\,\beta h \,r_n^2\,)^{1/\,2}\,.
\]
These inequalities are incompatible for $n$ large enough,
since $r_nt_n^{-1}\to +\infty$.

\end{proof}

\subsection{A closer look to the Lagrangian action estimate}

In the proof of our theorem we will make use of a slight modification
of the fundamental result stated in section \ref{subs: estimates}.
More precisely,
we refer to the existence of a curve joining two given configurations
in a given time,
for which an upper bound of its Lagrangian action is obtained.

The curve provided by the theorem,
a curve $\gamma\in\calC(x,y,\tau)$ is built in \cite{Ma1}
in canonical way from the data,
namely the configurations $x, y$ in $E^N$ and time $\tau>0$.
Moreover, for configurations $x=(r_1,\dots,r_N)$ and $y=(s_1,\dots,s_N)$
contained in a ball of radius $R$ of $E$, that is such that
\[
\norm{r_i-c}<R\;\;\;\textrm{ and }\;\;\;\norm{s_i-c}<R
\]
for some $c\in E$ and for all $i\in\calN$, the construction is such that
the configuration $\gamma(t)=(\gamma_1(t),\dots,\gamma_N(t))$
is contained in the concentric ball of radius $6NR$
for all $t\in [0,\tau]$.

As we will see, the canonical path $\gamma\in\calC(x,y,\tau)$ can
be constructed with the following additional condition:
if $x$ and $y$ are configurations with the same center of mass,
then the center of mass of the curve $\gamma(t)$ is constant
and equals the common center of the configurations $x$ and $y$.
However, to achieve this we must allow the configuration $\gamma(t)$
to be contained in a ball of radius twice as large, that is,
being able to assure that $\norm{\gamma_i(t)-c}<12NR$
for all $i\in\calN$ and for all $t\in [0,\tau]$.

One possible way to prove this,
which is not the one we will address,
consists of adapting the proof offered in \cite{Ma1}.
Indeed, the construction of the curve $\gamma$
is obtained there by reparamterizing the concatenation
of two linear paths with a common vertex in an intermediate configuration $z$.
To achieve the estimation of the action,
that intermediate configuration $z=(z_1,\dots,z_N)$ must be required
to have all mutual distances greater than or equal to $6R$.
In order for the argument to be valid even in the case dim E = 1,
a collinear $z$ configuration is used,
in which the body $z_1$ is in the center of the ball,
and the distance between $z_i$ and $z_{i+1}$ is exactly $ 6R $.
This explains the upper bound $6NR$ for the radius of the ball containing
the curve $\gamma$.
It is not difficult to verify that
the additional condition we are looking for,
is trivially obtained if we translate the intermediate configuration $z$
in such a way that its center of mass
coincides with the midpoint of the segment
joining the center of $x$ and the center of $y$.
It is also observed that the path obtained with this modification
could not be contained in the ball of radius $6NR$,
but is necessarily contained in the ball of radius $12NR$.

Instead,
to avoid these verifications for the reader,
we propose the use of the following simple lemma.

\begin{lemma}
\label{lemma: curva centrada}
Let $x,y\in E^N$, $\tau>0$ and $\gamma\in\calC(x,y,\tau)$.
If $\xi(t)\in E$ denotes the center of mass of
$\gamma(t)=(\gamma_1(t),\dots,\gamma_N(t))$
for each $t\in [0,\tau]$ then the curve $\gamma'\in\calC(x,y,\tau)$
defined by
\[
\gamma'_i(t)=
\frac{(\tau-t)}{\tau}\,\xi(0)+\frac{t}{\tau}\,\xi(\tau) \;+\;
\gamma_i(t)-\xi(t) 
\]
satisfies $\calA(\gamma')\leq\calA(\gamma)$ and moreover,
if $B\subset E$ is a ball containing $\gamma(t)$ for all $t\in [0,\tau]$
then $\gamma'(t)$ is all the time contained
in the ball concentric with $B$ and twice the radius.
\end{lemma}

\begin{proof}
Consider the trivial partition of $E^N$ with only one cluster $A=\calN$.
In this case, the orthogonal decomposition of $E^N$
is simply reduced to identifying each configuration $x=(r_1,\dots,r_N)$
with the pair $(y,(r_1-y,\dots,r_N-y))$,
being $y\in E$ the center of mass of the configuration $x$.

Using the corresponding orthogonal decomposition of $E^N$,
that is to say $E^N=E\times E^N_0$,
and writing $\sigma(t)=(\gamma_1(t)-\xi(t),\dots,\gamma_N(t)-\xi(t))$
for $t\in [0,\tau]$, we have that $\gamma(t)=(\xi(t),\sigma(t))$ and that
$\gamma'(t)=(\xi'(t),\sigma(t))$,
where
\[
\xi'(t)=\frac{(\tau-t)}{\tau}\,\xi(0)+\frac{t}{\tau}\,\xi(\tau)
\]
is the linear path from the center of $x$ to the center of $y$.

Therefore, using now the decomposition for the action, we get
\[
\calA(\gamma)-\calA(\gamma')=
\int_0^\tau\frac{1}{2}M\norms{\dot\xi(t)}^2dt\,-\,
\int_0^\tau\frac{1}{2}M\norms{\dot\xi'(t)}^2dt\,
\]
where $M$ is the total mass of the system.
We conclude that $\calA(\gamma)\geq\calA(\gamma')$
since the Bunyakowsky's inequality implies that
\[
\int_0^\tau\norms{\dot\xi(t)}^2dt\;\geq\;
\frac{1}{\tau}\,\left(\int_0^\tau\norms{\dot\xi(t)}dt\right)^2\;\geq\;
\frac{1}{\tau}\,\norms{\xi(\tau)-\xi(0)}^2
\]
with equality holding if and only if $\xi(t)=\xi'(t)$ for all $t\in [0,\tau]$.
Concerning the second assertion,
its proof is immediate since if $\gamma(t)$
is contained in some ball $B\subset E$ for all $t\in [0,\tau]$, then
both $\xi(t)$ and $\xi'(t)$ are in $B$ all the time.

\end{proof}

\subsection{Proof of Theorem \ref{mainthm0}}

\begin{proof}
Let $x:[0,+\infty)\to \Omega$ be a free time
minimizer of $\calA_h$, with $h\geq 0$, say $x(t)=(r_1(t),\dots,r_N(t))$.
By Lemma \ref{lemma: FTMnoSH},
we know that $x(t)$ is not superhyperbolic,
and that $x(t)=ta + O(t^{2/3})$.

In the case $h=0$,
we already know that $x(t)$ is a completely parabolic motion,
see \cite{DaLMa}.
Some times this is expressed by the condition $\dot x(t)\to 0$,
but equivalently we must have $r_{ij}(t)\approx t^{2/3}$ for all $i\neq j$.
It remains then to consider the case $h>0$, in which $a\neq 0$.
However, the following proof works in all cases, that is, even when $h=0$.

Let us write $a=(a_1,\dots,a_N)$, and call $\calP$
the natural partition associated to the motion $x$,
that is $i\sim j$ if and only if $a_i=a_j$.
According to the description of the asymptotic behaviour
of the clusters also given by Marchal and Saari and recalled
in section \ref{subsect:cluster-energy},
for each cluster $A\in\calP$ either the asymptotic energy
is $h_A=0$ and $r_{ij}(t)\approx t^{2/3}$ for any $i,j\in A$ with $i\neq j$,
or $h_A<0$ and the minimal distance
$r_A(t)=\min\set{r_{ij}(t)\mid i\neq j,\;i,j\in A}$ is bounded.
Therefore, it suffices to prove that $h_A=0$ for all $A\in\calP$ or
equivalently, that $r_A(t)$ is unbounded for all $A\in\calP$. 

The cluster decomposition associated to our partition will be denoted by
$d(x(t))=(y(t),(z_A(t))_{A\in\calP})$.
Recall that $y(t)=(y_A(t))_{A\in\calP}$ is the configuration in $E^\calP$
composed by the cluster centers,
and that $z_A(t)=(s_i(t))_{i\in A}\in\EA0$
is given for each cluster $A\in\calP$ by the relative positions
$s_i(t)=r_i(t)-y_A(t)$.

Now we use the fact that $x(t)=ta+O(t^{2/3})$.
This implies that there is $R>0$ such that
$\norm{r_i(t)-ta_i}<Rt^{2/3}$ for all $t>t_0$ for some $t_0>0$.
Thus for each cluster $A$ we also have $\norm{y_A(t)-ta_i}<Rt^{2/3}$
whenever $i\in A$,
hence $\norm{s_i(t)}<2Rt^{2/3}$ for all $i\in\calN$ and every $t>t_0$.
This means that the configuration $z_A(t)\in\EA0$ is contained in the ball
with center at the origin of $E$ and radius $2Rt^{2/3}$ for $t>t_0$.

Suppose that there is a cluster $A\in\calP$ and $\rho>0$ such that
$r_{A}(t)<\rho$ for all $t\geq 0$.
Then we deduce that 
\[
U_A(x(t))\geq k=\frac{m_0^2}{\rho}>0
\]
for all $t\geq 0$, where $m_0=\min\set{m_1,\dots,m_N}$.
In particular we can state that,
\[
\int_t^{2t}L_A(z_A(s),\dot z_A(s))\,ds\geq kt
\]
for any $t>0$, where $L_A:T\EA0\to\R$ is the Lagrangian
given by
\[
L_A(z,w)=\frac{1}{2}\norm{w}^2\,+\,U_A(z)\,.
\]

We will show that for $t>t_0$ big enough,
this lower bound is incompatible
with the minimizing property of the motion $x(t)$.
To do this, for each $t>t_0$ we consider a new curve
$\gamma_t\in\calC(x(t),x(2t),t)$ as follows.
For the cluster $A\in\calP$ that is supposed to be non expanding,
we choose a curve $\sigma_A(s)$ in $\EA0$,
joining $z_A(t)$ with $z_A(2t)$ and defined for $s\in [t,2t]$,
in such a way that we have
\[
\int_t^{2t}L_A(\sigma_A(s),\dot \sigma_A(s))\,ds\leq
\alpha_0\,m_A\,n_A^4\,\frac{R_{2t}^{\,2}}{t}\,+\,
\beta_0\,m_A^2\,n_A^4\,\frac{t}{R_{2t}}\,,
\]
where $m_A$ is the total mass of the cluster,
$n_A=\#A$ is the number of bodies in the cluster,
$R_t=2R\,t^{2/3}$ is the radius of the ball containing
both configurations $z_A(t)$ and $z_A(2t)$ in $\EA0$,
and $\alpha_0$ and $\beta_0$ are the positive constants
given by the fundamental theorem stated in section \ref{subs: estimates}.
Note that here we are applying this theorem
to the Lagrangian $L_A$ in $E^A$.
Moreover, according to Lemma \ref{lemma: curva centrada},
we can assume that $\sigma_A(s)\in\EA0$ for all $s\in [t,2t]$,
and that it is contained in the ball
with center at the origin and radius $12\,n_AR_{2t}$.
Thus we have
\[
\int_t^{2t}L_A(\sigma_A(s),\dot \sigma_A(s)\,ds\leq
k'\,t^{1/3}
\]
for some constant $k'>0$, and this for all $t>t_0$.
Then we define the curve $\gamma_t\in\calC(x(t),x(2t),t)$
by replacing in the cluster decomposition of the original curve
the component $z_A(s)$ by the curve $\sigma_A(s)$.

In order to get a contradiction, finally we compare de Lagrangian action
of the curve $x\mid_{[t,2t]}$ with that of the new curve $\gamma_t$.
Using the decomposition of the action we see,
after cancellation of the common terms, that
\[
\calA(\gamma_t)-\calA(x\mid_{[t,2t]})\;\leq\;
k'\,t^{1/3}\,-\,k\,t \;+\;
\int_t^{2t}W_\calP(\gamma_t(s))-W_\calP(x(s))\,ds\,.
\]

Let us fix for a while the value of $t>0$.
If we write $\gamma_t(s)=(r'_1(s),\dots,r'_N(s))$ then we have that
\[
W_\calP(\gamma_t(s))-W_\calP(x(s))=
\sum_{i\in A} \sum_{j\notin A}m_im_j\,
\left(
\frac{1}{\norm{r_i(s)-r_j(s)}}-\frac{1}{\norm{r'_i(s)-r_j(s)}}
\right).
\]
By the previous considerations, for any $i\in A$ and $j\notin A$ we also have
\begin{eqnarray*}
\norm{r_i(s)-r_j(s)} & \geq &
\norm{sa_i-sa_j}-\norm{r_i(s)-sa_i}-\norm{r_j(s)-sa_j}\\
&\geq & s\norm{a_{ij}}-2R\,s^{2/3}
\end{eqnarray*}
as well as
\begin{eqnarray*}
\norm{r'_i(s)-r_j(s)} & \geq &
\norm{sa_i-sa_j}-\norm{r'_i(s)-y_A(s)}\\
&&\;\;-\norm{y_A(s)-sa_i}-\norm{r_j(s)-sa_j}\\
&\geq & s\norm{a_{ij}}-12\,n_A\,R\,s^{2/3}- 2R\,s^{2/3}\,.
\end{eqnarray*}
Therefore we deduce that there is a constant $k''>0$ and $t_1>t_0$
such that
\[
\abs{W_\calP(\gamma_t(s))-W_\calP(x(s))}\;\leq\;\frac{\;k''}{s}
\]
for all $s\in [t,2t]$, whenever $t\geq t_1$, which in turn implies that
\[
\calA(\gamma_t)-\calA(x\mid_{[t,2t]})\;\leq\;
k'\,t^{1/3}\,-\,k\,t \,+\,k''\log(2)\,.
\]
For $t>t_1$ big enough we have that this upper bound becomes negative,
meaning that $\calA(x\mid_{[t,2t]})>\phi(x(t),x(2t),t)$.
But clearly this is in contradiction with the hypothesis that $x$ is a free time
minimizer of $\calA_h$.

\end{proof}

\section*{Acknowledgments}

The first author is grateful to Cátedras Conacyt program of
Consejo Nacional de Ciencia y Tecnología del Gobierno de México.
The second author thanks Grupo CSIC 618 UdelaR,
Uruguay, as well as the mathematics department of
Cinvestav (Instituto Politécnico Nacional, México)
 for the hospitality during part of the development of this work.

%%%%%%%%%%%%%%%%%%%%%%%%%%%%%%%%%%%%%%%%%%
%%%%%%%%%%%%%%%%%%%%%%%%%%%%%%%%%%%%%%%%%%
%%%%%%%%%%%%%%%%%%%%%%%%%%%%%%%%%%%%%%%%%%
%%%%%%%%%%%%%%%%%%%%%%%%%%%%%%%%%%%%%%%%%%
%%%%%%%%%%%%%%%%%%%%%%%%%%%%%%%%%%%%%%%%%%


\begin{thebibliography}{Coh93}

\bibitem{BarSec}
  \textsc{V.~Barutello, S.~Secchi.}
  Morse index properties of colliding solutions to the
  $N$-body problem,
  \emph{Ann. Inst. H. Poincaré Anal. Non Linéaire}
  \textbf{25} n.3 (2008), 539--565.

\bibitem{Bur}
  \textsc{J.~M.~Burgos},
  Existence of partially hyperbolic motions in the
  $N$-body problem,
  \emph{preprint in arXiv}.

\bibitem{Cha1}
  \textsc{J.~Chazy},
  Sur certaines trajectoires du problème des n corps,
  \emph{Bull. astronom.}
  \textbf{35} (1918), 321--389.

\bibitem{Cha2}
  \textsc{J.~Chazy},
  Sur l'allure du mouvement dans le problème des trois corps
  quand le temps croît indéfiniment, 
  \emph{Ann. Sci. E.N.S. (3-ème série)}
  \textbf{39} (1922), 29--130.

\bibitem{Che}
  \textsc{A.~Chenciner},
  Action minimizing solutions of the Newtonian
  $n$-body problem: from homology to symmetry,
  \emph{Proceedings of the
  International Congress of Mathematicians Beijing 2002},
  vol. III, 279--294, Higher Ed. Press, Beijing, 2002.
  
\bibitem{CoPa}
  \textsc{G.~Contreras, G.~P.~Paternain},
  Connecting orbits between static classes for generic Lagrangian systems,
  \emph{Topology}
  \textbf{41} (2002), 645--666.
  
\bibitem{DaLMa}
  \textsc{A.~Da Luz, E.~Maderna},
  On the free time minimizers of the Newtonian $N$-body problem,
  \emph{Math. Proc. Camb. Phil. Soc.}
  \textbf{156} (2014), 209--227.
  
\bibitem{DMMY}
  \textsc{N.~Duignan, R.~Moeckel,
  R.~Montgomery, G.~Yu},
  Chazy-type asymptotics and hyperbolic scattering
  for the $n$-body problem,
  \emph{Arch. Ration. Mech. Anal.}
  \textbf{238} n.1 (2020), 255--297.

\bibitem{Fat}
	\textsc{A.~Fathi},
	Weak KAM Theory:
	the connection between Aubry-Mather theory and
	viscosity solutions of the Hamilton-Jacobi equation,
	\emph{Proceedings of the International Congress of
	Mathematicians Seoul 2014}, vol. III, 597--621.

\bibitem{FerTer}
  \textsc{D.~Ferrario, S.~Terracini},
  On the existence of collisionless equivariant minimizers
  for the classical $n$-body problem,
  \emph{Invent. Math.}
  \textbf{155} n.2 (2004), 305--362.

\bibitem{Ma1}
  \textsc{E.~Maderna},
  On weak KAM theory for $N$-body problems,
  \emph{Ergod. Th. Dynam. Sys.}
  \textbf{32} n.3 (2012), 1019--1041.

\bibitem{Ma2}
  \textsc{E.~Maderna},
  Minimizing configurations and Hamilton-Jacobi equations
  of homogeneous $N$-body problems,
  \emph{Regul. Chaotic Dyn.}
  \textbf{18} n.6 (2013), 656--673.

%\bibitem{MaVe1}
% \textsc{E.~Maderna, A.~Venturelli},
% Globally Minimizing Parabolic Motions
%  in the Newtonian N -body Problem,
% \emph{Arch. Ration. Mech. Anal.}
% \textbf{194} n.1 (2009), 283--313.

\bibitem{MaVe2}
  \textsc{E.~Maderna, A.~Venturelli},
  Viscosity solutions and hyperbolic motions:
  A new PDE method for the $N$-body problem,
  \emph{Ann. of Math. (2)}
  \textbf{192} n.2 (2020), 499--550.
  
\bibitem{Man}
  \textsc{R.~Mañé},
  Lagrangian flows: the dynamics of globally minimizing orbits,
  \emph{Bol. Soc. Bras. Mat.}
  \textbf{28} n.2 (1997), 141--153.

\bibitem{Mar}
  \textsc{C.~Marchal},
  How the method of minimization of action avoids singularities,
  \emph{Celestial Mech. Dynam. Astronom.}
  \textbf{83} (2002), 325--353.

\bibitem{MarSaa}
  \textsc{C.~Marchal, D.~Saari},
  On the final evolution of the $n$-body problem,
  \emph{J. Differential Equations}
  \textbf{20} n.1 (1976), 150--186.

\bibitem{McG}
	\textsc{R.~McGehee},
	Von Zeipel's theorem on singularities in celestial mecanics,
	\emph{Expo. Math.}
	\textbf{4} (1986), 335--345.
	
\bibitem{MoMoSa}
  \textsc{R.~Moeckel, R.~Montgomery, H.~Sánchez Morgado},
  Free time minimizers for the three-body problem,
  \emph{Clest. Mech. Dyn. Astr.}
  \textbf{130} n.3 (2018).

\bibitem{PerSan}
  \textsc{B.~Percino-Figueroa} and \textsc{H.~Sánchez-Morgado},
  Busemann functions for the $N$-body problem,
  \emph{Arch. Ration. Mech. Anal.}
  \textbf{213} no.3 (2014), 981--991.

\bibitem{Xue}
  \textsc{J.~Xue},
  Non-collision singularities in a planar $4$-body problem,
  \emph{Acta Mathematica}
  \textbf{224} n.2 (2020), 253--388.

\bibitem{Zei}
  \textsc{H.~von Zeipel},
  Sur les singularités du problème des n corps,
  \emph{Ark. Math. Astr. Fys.}, n.4 (1908), 1--4.  

\end{thebibliography}
\end{document}